\documentclass{amsart}

\usepackage{amsmath}
\usepackage{amssymb}  
\usepackage{amsthm}   
\usepackage{url}
\usepackage{color}

\DeclareMathOperator{\Aut}{Aut}

\newcommand{\tth}{^{\operatorname{th}}}

\theoremstyle{plain}
\newtheorem{thm}{Theorem}
\theoremstyle{definition}
\newtheorem{defn}[thm]{Definition}
\theoremstyle{remark}
\newtheorem*{rem}{Remark}

\begin{document}
    \title{Finding Rational Periodic Points on Wehler K3 Surfaces}
    \author[Hutz]{Benjamin Hutz}
    \address{Department of Mathematics and Computer Science,
            Amherst College,
            Amherst, MA}
    \email{bhutz@amherst.edu}

    \thanks{The author would like to thank his advisor Joseph Silverman for his many insightful suggestions and ideas.}




    \begin{abstract}
        This article examines dynamical systems on a class of K3 surfaces in $\mathbb{P}^{2} \times \mathbb{P}^{2}$ with an infinite automorphism group.  In particular, this article develops an algorithm to find $\mathbb{Q}$-rational periodic points using information modulo $p$ for various primes $p$.  The algorithm is applied to exhibit K3 surfaces with $\mathbb{Q}$-rational periodic points of primitive period $1,\ldots,16$.  A portion of the algorithm is then used to determine the Riemann zeta function modulo $3$ of a particular K3 surface and find a family of K3 surfaces with Picard number two.
    \end{abstract}

\maketitle

     \section{Introduction}
        Let $\phi: V \subseteq \mathbb{P}^N_K \to V$ be a morphism on a variety $V$.  We denote the $n\tth$ iterate of $\phi$ as
        \begin{equation*}
            \phi^n = \phi \circ \phi^{n-1}.
        \end{equation*}
        We say that $P \in V$ is \emph{preperiodic} if $\{\phi^i(P) \mid i \in \mathbb{N}\}$ is finite.  In particular. we say that $P$ is a \emph{periodic point of period $n$ for $\phi$} if $\phi^n(P) = P$
        and of \emph{primitive period $n$} if, in addition, for all $m < n$, $\phi^m(P) \neq P$.

        By work of Northcott \cite{Northcott}, there are only finitely many rational periodic points for a morphism of projective space.  In particular, this implies that the possible primitive period of a rational periodic point is bounded.  Morton and Silverman \cite{Silverman7} conjecture for morphisms of projective space that this bound should depend only the degree of the number field, the dimension of the space, and the degree of the morphism.  Fakhruddin \cite{Fakhruddin} has shown that this conjecture implies a similar statement about morphisms of certain algebraic varieties.  This article examines dynamical systems on a class of surfaces in $\mathbb{P}^{2} \times \mathbb{P}^{2}$ with an infinite automorphism group studied by Wehler \cite{Wehler}.  Silverman \cite[Theorem 1.1]{Silverman5} shows that the set of rational (pre)periodic points for the composition of the two non-commuting involutions on these K3 surfaces is a set of bounded height and hence finite.  Call and Silverman \cite[Appendix]{Silverman4} give explicit formulas for computing the two non-commuting involutions on these K3 surfaces.  This article addresses the following topics:
        \begin{enumerate}
            \item find $\mathbb{Q}$-rational periodic points on Wehler K3 surfaces using information modulo $p$ for various primes $p \in \mathbb{Q}$.
            \item determine the number of rational points on a Wehler K3 surface in the field with $p^m$ elements, denoted $\mathbb{F}_{p^m}$.
        \end{enumerate}
        In the course of addressing these topics, we produce an explicit family of K3 surfaces with Picard number two.  Note that van Luijk \cite{vanluijk} constructs explicit K3 surfaces with Picard number one and infinitely many rational points using more general methods than discussed in this article.  However, since his examples have an elliptic fibration his methods are not wholly applicable to Wehler K3 surfaces.  Additionally, the orbits of curves and points have been studied on these K3 surfaces by Baragar in \cite{Baragar3,Baragar2}.

        This article proceeds as follows.  In Section \ref{sect2}, we define and recall some basic information about Wehler's class of K3 surfaces.  In Section \ref{sect3}, we address finding $\mathbb{Q}$-rational periodic points from data modulo primes.  The key fact used to find rational periodic points is that if there is a rational periodic point, then the primitive period of its reduction modulo a prime of good reduction must divide the primitive period of the rational point.  While this method may return false positives due to checking only a finite number of primes, it provides an effective way to search for rational periodic points.  This algorithm is used successfully in Pari/gp to find surfaces with $\mathbb{Q}$-rational points of primitive period $1,\ldots, 16$.  In Section \ref{sect6}, we use the number of $\mathbb{F}_{p^m}$ rational points with $m=1,\ldots, 11$ for a particular Wehler K3 surface to determine the Riemann zeta function of the surface using the Weil conjectures.  An example is computed for $p=3$ and shown to have Picard number two which produces an infinite family of Wehler K3 surfaces with Picard number two.  Much of this work is from the author's doctoral thesis \cite[Chapters 6 and 7]{Hutz}.

        There are also several similarly constructed classes of surfaces in Baragar \cite{Baragar}, Billard \cite{Billard}, and McMullen \cite{McMullen} that may allow a similar type of explicit investigation.

        A PARI/gp script to search for a $\mathbb{Q}$-rational primitive periodic point on a Wehler K3 surface and to check for the possible primitive periods of a periodic point on a Wehler K3 surface is available on the author's homepage.

     \section{Wehler's K3 surfaces} \label{sect2}
        A Wehler K3 surface $S \subset \mathbb{P}^{2}_x \times \mathbb{P}^{2}_y$ is
        a smooth surface given by the intersection of an effective divisor of degree (1,1) and an effective
        divisor of degree (2,2).  In other words,  let $([x_0,x_1,x_2],[y_0,y_1,y_2]) = (\textbf{x},\textbf{y})$ be the coordinates for $\mathbb{P}^2_x \times \mathbb{P}^2_y$, then
        $S$ is the locus described by $L(\textbf{x},\textbf{y})=Q(\textbf{x},\textbf{y})=0$ for
        \begin{equation*}
            L(\textbf{x},\textbf{y}) = \sum_{0\leq i,j \leq 2} a_{ij}x_iy_j \qquad
            Q(\textbf{x},\textbf{y}) = \sum_{0\leq i,j,k,l \leq 2} b_{ijkl}x_ix_jy_k y_l.
        \end{equation*}
        \begin{thm}\label{thm1} \textup{\cite[Theorem 2.9]{Wehler}}
            A general K3 surface formed as the vanishing locus of a degree $(1,1)$ and a degree $(2,2)$ effective divisor has Picard number two and an infinite automorphism group.
        \end{thm}
        When working with a particular surface, we will specify the vectors
        \begin{align*}
            L &= [a_{ij}] = [a_{00},a_{01},a_{02},a_{10},\ldots, a_{22}] \quad \text{and} \\
            Q &= [b_{ijkl}] = [b_{0000},b_{0001},b_{0002},b_{0010},\ldots,b_{2222}].
        \end{align*}
        The natural projections
        \begin{equation*}
            \rho_x: \mathbb{P}^{2}_x \times \mathbb{P}^2_y \to \mathbb{P}^2_x, \quad
            \rho_y: \mathbb{P}^{2}_x \times \mathbb{P}^2_y \to \mathbb{P}^2_y
        \end{equation*}
        induce two projection maps
        \begin{equation*}
            p_x: S \to \mathbb{P}^2_x, \quad
            p_y: S \to \mathbb{P}^2_y.
        \end{equation*}
        The projections $p_x$ and $p_y$ are in general double covers, allowing us to define two involutions of $S$, $\sigma_x$ and $\sigma_y$, respectively.  To define $\sigma_x$, we consider $p^{-1}_x(\textbf{x}) = \{(\textbf{x},\textbf{y}),(\textbf{x},\textbf{y}')\}$ and define $\sigma_x((\textbf{x},\textbf{y})) = (\textbf{x},\textbf{y}')$.  To define $\sigma_y$, we consider $p^{-1}_y(\textbf{y}) = \{(\textbf{x},\textbf{y}),(\textbf{x}',\textbf{y})\}$ and define $\sigma_y((\textbf{x},\textbf{y})) = (\textbf{x}',\textbf{y})$.  The maps $\sigma_x$ and $\sigma_y$ are in general just rational maps.  However, if $S=V(L,Q)$ is smooth, Call and Silverman \cite[Proposition 1.2]{Silverman4} show that $\sigma_x$ and $\sigma_y$ are morphisms of $S$ if and only if $S$ has no degenerate fibers, fibers of positive dimension.  We call a surface with no degenerate fibers a \emph{non-degenerate} surface.
        \begin{rem}
            In the case where $S$ is smooth and degenerate, Baragar \cite[Section 2]{Baragar2} shows how to extend $\sigma_x$ and $\sigma_y$ to morphisms defined also on the degenerate fiber.  These extended morphisms are not treated here, but offer an interesting area of future study.
        \end{rem}
        We adopt the notation from \cite{Silverman4} and define
        \begin{align*}
            L_j^{x} &= \text{ the coefficient of } y_j \text{ in } L(\textbf{x},\textbf{y}), \\
            L_j^{y} &= \text{ the coefficient of } x_j \text{ in } L(\textbf{x},\textbf{y}),\\
            Q_{kl}^{x} &= \text{ the coefficient of } y_ky_l \text{ in } Q(\textbf{x},\textbf{y}),\\
            Q_{ij}^{y} &= \text{ the coefficient of } x_ix_j \text{ in } Q(\textbf{x},\textbf{y}),\\
            G_k^{\ast} &= (L_{j}^{\ast})^2Q_{ii}^{\ast} - L_i^{\ast}L_j^{\ast}Q_{ij}^{\ast} + (L_i^{\ast})^2Q_{jj}^{\ast},\\
            H_{ij}^{\ast} &= 2L_i^{\ast}L_j^{\ast}Q_{kk}^{\ast} - L_i^{\ast}L_k^{\ast}Q_{jk}^{\ast} - L_j^{\ast}L_k^{\ast}Q_{ik}^{\ast} + (L_k^{\ast})^2Q_{ij}^{\ast}.
        \end{align*}
        for $(i,j,k)$ some permutation of the indices $\{0,1,2\}$ and $\ast$ replaced by either $x$ or $y$.

     \section{Periodic Points}\label{sect3}
        We use good reduction information to develop an algorithm to find Wehler K3 surfaces with $\mathbb{Q}$-rational periodic points of specified primitive period.
        \begin{thm}
            Let $K$ be a number field.  Assume that
            $S \subset \mathbb{P}_x^{2}\times \mathbb{P}_y^{2}$ is a Wehler K3 surface defined over $K$.
            Let $\phi=\sigma_x\circ \sigma_y \in \Aut(S)$, let $P \in S(K)$ be a periodic point for $\phi$, and let $\mathfrak{P} \subset K$ be the set of primes where $\phi$ has good reduction.  Define
            \begin{enumerate}
                \item[] $n =$ primitive period of $P$ for $\phi$, and
                \item[] $m_{\mathfrak{p}} =$ primitive period of $\overline{P}$ for $\overline{\phi}$ for $\mathfrak{p} \in \mathfrak{P}$.
            \end{enumerate}

            Then $m_{\mathfrak{p}} \mid n$ for all primes $\mathfrak{p} \in \mathfrak{P}$ and $n=m_{\mathfrak{p}}$ for all but finitely many primes $\mathfrak{p} \in \mathfrak{P}$.
        \end{thm}
        \begin{proof}
            To show divisibility we know for primes of good reduction
            \begin{equation*}
                \overline{\phi^n(P)} = \overline{\phi}^n(\overline{P}) = \overline{P}
            \end{equation*}
            so we have $m_{\mathfrak{p}} \leq n$.  Since the smallest period is the greatest common divisor of all of the periods, we must have $m_{\mathfrak{p}} \mid n$.

            For the second statement, let $x_i(P)$ and $y_i(P)$ for $0 \leq i \leq 2$ be the coordinates of $P$.  The assumption that $\phi^n(P) =P$ is equivalent to the system of equations
            \begin{align*}
                c_1x_0(\phi^n(P)) &- x_0(P) =0, \quad c_1x_1(\phi^n(P)) - x_1(P) =0, \quad c_1x_2(\phi^n(P)) - x_2(P) =0,\\
                c_2y_0(\phi^n(P)) &- y_0(P) =0, \quad c_2y_1(\phi^n(P)) - y_1(P) =0, \quad \text{and} \quad c_2y_2(\phi^n(P)) - y_2(P) =0,
            \end{align*}
            for some nonzero constants $c_1,c_2 \in K$.
            We know that there are only finitely many integers $1 \leq d< n$ and each value $c_1x_i(\phi^d(P))-x_i(P)$ and $c_2y_i(\phi^d(P))-y_i(P)$ for $0 \leq i \leq 2$  can be factored into finitely many primes in $K$, so there can only be finitely many primes such that $\phi^d(P) =P$ for $d < n$.
        \end{proof}
        The algorithm proceeds as follows to find a $\mathbb{Q}$-rational periodic point with specified primitive period $n$ on a given surface $S$:
        \begin{enumerate}
            \item Chose a prime $p \in \mathbb{Q}$.
            \item Find all of the points in $S(\mathbb{F}_p)$. (Discussed in Section \ref{ssect2})
            \item Find all of the cycles in $S(\mathbb{F}_p)$.
                \subitem Exit if $S$ is degenerate.
            \item Save all points with $m \mid n$ for $p$.
                \subitem Exit if there are none.
            \item Repeat for the next prime until the specified number of primes is completed.
            \item Use the Chinese Remainder Theorem on the coordinates of the
                saved points to find $P \in S(\mathbb{Q})$.
        \end{enumerate}
        There are two points worth further discussion: What if $S=V(L,Q)$ is singular and thus not a K3 surface? Which points do you use for the CRT?

        First we address singularity.  Since we are computing the forward image of every point in $S(\mathbb{F}_p)$ we can easily check for degenerate fibers since by \cite{Silverman4} a fiber $S_a^{\ast} = p_{\ast}^{-1}(\textbf{a})$ is degenerate if and only if
        \begin{equation*}
            G_0^{\ast}(\textbf{a}) = G_1^{\ast}(\textbf{a}) = G_2^{\ast}(\textbf{a})
            = H_{01}^{\ast}(\textbf{a}) = H_{02}^{\ast}(\textbf{a}) = H_{12}^{\ast}(\textbf{a}) = 0.
        \end{equation*}
        From \cite[Proposition 2.5]{Silverman4} this leaves the only possible singular points as fixed points of $\sigma_x$ and $\sigma_y$.  Therefore, if we are searching for a periodic point with primitive period $n>1$ the singular fixed points have minimal effect on the algorithm.  So we content ourselves to only check for degenerate fibers.  When a surface $S=V(L,Q)$ is found to have a periodic point we can then check to see if it is nonsingular.

        Determining which of the saved points to use in the Chinese Remainder Theorem is a more practical issue.  We could try every possible combination of points with $m \mid n$ modulo each prime $p$, however this may not be the most efficient strategy.  Since we are only interested in finding a periodic point and not necessarily all periodic points and are willing to accept missing some periodic points for a faster algorithm, we choose only points from primes where there is a single point with $m=n$ or primes where there is only a single point with $m \mid n$.

        The above algorithm was used to find Wehler K3 surfaces with $\mathbb{Q}$-rational periodic points of primitive period $1,\ldots,16$.  Table \ref{table2} and Table \ref{table3} give examples of surfaces with a point of each primitive period.  Note that the surface with a $6$-periodic point also has a $2$-periodic point.

\begin{table}
    \caption{Primitive Periods 1-8}
    \label{table2}
    \begin{center}
    \begin{tabular}{|c|l|}
        \hline
        $n$& Point $P$ and surface $S = V(L,Q)$ \\
        \hline
         $1$ & $P=[[0,0,1],[1,0,1]]$\\
          & $L:[-1, 1, 1, 0, -1, -2, 0, 1, 0]$\\
          & $Q:[-2, -1, -1, 1, -2, 0, 1, -2, -1, 0, 1, -1, 1, 1, 1, 0, -2, -2, -1, -2, -1, -2, -2, 1,$\\
          & \qquad $0, -2,-1, 1, -1, -1, 1, 1, -2, 1, 1, 1]$\\
        \hline
        $2$ & $P = [[3,1,3],[1,0,0]]$ \\
          & $L:[1, -3, 2, 0, 1, -3, -1, -3, -1]$\\
          & $Q:[-2,-4,-2,-2,-3,2,0,-2,0,1,-2,0,1,2,1,-4,0,0,3,2,-3,1,-2,1,-4,-2,$\\
          & \qquad $1,1,2,2,2,-1,1,2,-4,0]$ \\
        \hline
        $3$ & $P = [[1,-1,0],[1,0,0]]$ \\
          & $L:[0, -2, 0, 0, -2, 1, 1, -1, -2]$ \\
          & $Q:[0, 1, 0, 1, -1, -2, -1, 0, 0, 0, -1, -2,-2, 1, -1, -1, -1, 1,-1, -2, -1, -1,0, 1,$\\
          &  \qquad $1, 0, -2,-2, 1,1, -2, -2, -2, -2, 1, -2]$ \\
        \hline
        $4$ & $P = [[1, 0, -1], [1, -1, 0]]$ \\
          & $L:[-1, -1, 0, -1, 0, -1, 0, 0, -1]$\\
          & $Q:[-1, 0, -1, 0, 0, -1, -1, 0, -1, 0, 0, -1,-1,0, 0, 0, 0, 0, 0, 0, 0, -1, 0, -1, 0, -1,$\\
          & \qquad $0, 0, -1, 0, 0, -1, 0, -1, 0, 0]$ \\
        \hline
        $5$ & $P = [[1,-1,0],[1,0,1]]$ \\
          & $L:[0, 0, -1, -1, 0, 0, 0, -1, -1]$\\
          & $Q:[-1, -1, 0, 0, 0, 0, -1, 0, 0, 0,-1, -1, 0, 0, -1, 0, 0, 0, 0, -1, 0, -1, 0, -1,-1,0,$\\
          & \qquad $-1, 0, -1, -1, -1, -1, 0, 0, 0, -1]$ \\
        \hline
        $6$ & $P_6 =[[1,0,0],[1,0,0]] \quad P_2 = [[0,1,0],[1,0,-1]]$\\
          & $L:[0, 0, -1, -1, 0, -1, 0, -1, 0]$\\
          & $Q:[0, 0, 0, -1, 0, -1, -1, 0, 0, 0, -1, -1, -1,-1, 0, 0, -1, -1, 0, -1, 0, 0, 0, 0, 0,-1,$\\
          & \qquad $0, 0,-1, 0, 0, 0, -1, 0, 0, 0]$\\
        \hline
        $7$ & $P = [[1,-2,1],[1,-1,-1]]$\\
          & $L:[0, -1, 0, -1, 0, -1, -1, 0, 0]$\\
          & $Q:[-1, -1, -1, -1, 0, 0, 0, 0, -1, -1, -1, -1,-1, 0, 0, 0, 0, 0, 0, 0, -1, -1, 0, 0,0,$\\
          & \qquad $-1, -1, -1,0, -1, 0, 0, 0, -1, -1, -1]$ \\
        \hline
          $8$ & $P= [[1,0,0],[1,0,0]]$\\
           & $L:[0, 0, -1, -1, -1, 0, 0, -1, -1]$\\
           & $Q:[0, -1, 0, -1, -1, 0, 0, 0, 0, -1, -1, -1, -1, -1, 0, 0, -1, -1, 0, 0, -1, -1,-1, 0,$\\
           & \qquad $0,0, -1, -1, -1, 0, 0, -1, -1, 0, 0,-1]$\\
        \hline
    \end{tabular}
    \end{center}
\end{table}

\begin{table}
    \caption{Primitive Periods 9-16}
    \label{table3}
    \begin{center}
    \begin{tabular}{|c|l|}
        \hline
        $n$& Point $P$ and surface $S = V(L,Q)$ \\
        \hline
        $9$ & $P = [[0,1,0],[1,0,0]]$\\
          & $L:[-1, 0, 0, 0, 0, -1, 0, -1, 0]$\\
          & $Q:[-1, 0, 0, -1, -1, -1, -1, 0, -1, -1, -1, -1,0, 0, -1, -1, -1, -1, 0, -1, -1, 0,$\\
          & \qquad $0, 0,-1,-1,0, -1,-1, 0, 0, -1, -1, -1, 0, 0]$\\
        \hline
        $10$ & $P = [[1,-1,0],[0,1,-1]]$\\
           & $L:[0, -1, 0, -1, -1, 0, 0, -1, -1]$\\
           & $Q:[-1, -1, 0, 0, 0, 0, 0, 0, 0, 0, -1, -1,-1, 0, 0,-1, -1, -1, 0, 0, 0, -1, -1, 0, -1,$\\
           & \qquad $0, -1, 0,-1, 0, -1, -1, -1, 0, 0, 0]$\\
        \hline
        $11$ & $P = [[1, 0, -1], [-2, 1, 0]]$\\
           & $L:[-1, 0, 0, -1, -1, -1, -1, 0, -1]$\\
           & $Q:[0, -1, 0, -1, -1, 0, 0, 0, -1, -1, 0, -1,0, 0, 0, 0, -1,0, 0, 0, 0, -1, -1, 0,-1,-1,$\\
           & \qquad $-1, 0, 0, -1, 0, 0, -1, -1, -1, -1]$ \\
        \hline
        $12$ & $P=[[2,1,-2], [0, -2, 1]]$ \\
           & $L:[-1, -1, -1, -1, -1, 0, 0, -1, 0]$ \\
           & $Q:[0, -1, 0, 0, 0, 0, -1, 0, 0, -1, 0, 0, 0, 0,-1, -1, -1, -1, -1, -1, 0, -1, -1, 0,$\\
           & \qquad $-1, -1,0,-1,0, -1, -1, 0, -1, -1, -1, -1]$ \\
        \hline
        $13$ & $P =  [[0,1,0],[1,-1,-1]]$\\
           & $L:[0, -1, 0, -1, -1, 0, -1, 0, -1]$\\
           & $Q:[-1,0,-1,-1,0,0,-1,0,-1,0,0,-1,0,-1,0,0,-1,-1,-1,0,-1,0,0,0,0,$\\
           & \qquad $-1,0,-1,-1,-1,0,-1,0,0,0,-1]$\\
        \hline
        $14$ & $P = [[0, -2, -1], [-1, -3, 2]]$\\
           & $L:[0, -1, -1, -1, 0, 0, 0, 0, -1]$\\
           & $Q:[0, 0, -1, -1, 0, 0, 0, -1, 0, 0, 0, 0, -1, -1,-1, -1, -1, 0, -1, -1, -1, 0, -1,$\\
           &\qquad $0, 0, -1, 0, 0, 0, -1, -1, 0, -1, -1, -1, 0]$\\
        \hline
        $15$ & $P = [[1, 1, 0], [1, -1, -1]]$\\
           & $L:[-1, -1, 0, -1, 0, -1, 0, -1, -1]$\\
           & $Q:[0, -1, -1, 0, 0, 0, -1, 0, 0, 0, -1, 0, 0, -1, 0, 0,-1, -1, -1, -1, -1, 0, 0, -1, 0,$\\
           & \qquad $0, -1, -1,0, 0, -1, -1, 0, 0, -1, 0]$\\
        \hline
        $16$ & $P = [[1, 0, -1], [1, -1, 0]]$\\
           & $L:[-1, -1, 0, -1, 0, 0, 0, 0, -1]$ \\
           & $Q:[0, 0, -1, -1, 0, 0, -1, 0, -1, 0, -1, -1, -1,0, 0, 0, -1, 0, -1, 0, -1, 0, 0, -1, 0,$\\
           & \qquad $0, -1,-1,0, 0, -1, -1, -1, 0, 0, 0]$\\
        \hline
    \end{tabular}
    \end{center}
\end{table}

    Upon examining these tables, the reader will quickly notice that nearly all of the surfaces have coefficients contained completely in $\{-1,0,1\}$.  Upon examining data from the first 30 primes for randomly generated $S=V(L,Q)$, surfaces with larger coefficients had comparatively few periodic points: $48,871$ out of $49,200$ surfaces ($99\%$) with coefficients in $[-128,128]$ did not have a periodic point of any order.  The remaining $329$ surfaces may or may not have a periodic point but would require data from more primes to determine the existence of a periodic point.  Whereas for surfaces with coefficients in $\{-1,0,1\}$, $7425$ out of $19424$ surfaces ($37\%$) did not have periodic points of any order, with the remaining $11999$ surfaces requiring more data to determine the existence of periodic points.

    \subsection{Finding $p_{x}^{-1}(P)$} \label{ssect2}
        The algorithm requires finding $S(\mathbb{F}_p)$ efficiently.  We do so by finding the $0$, $1$, or $2$ points in $p_x^{-1}(\textbf{x})$ for $\textbf{x} \in \mathbb{P}^{2}(\mathbb{F}_p)$.  We know from Call
        and Silverman \cite[Corollary 1.5]{Silverman4} that the two solutions
        $\{\textbf{y},\textbf{y}^{\prime}\}$ satisfy the following equations:
        \begin{align*}
            y_0y_0^{\prime} &= G^x_0(\textbf{x}) \\
            y_1y_1^{\prime} &= G^x_1(\textbf{x}) \\
            y_2y_2^{\prime} &= G^x_2(\textbf{x}) \\
            y_0y_1^{\prime} + y_0^{\prime}y_1 &= -H^x_{01}(\textbf{x}) \\
            y_0y_2^{\prime} + y_0^{\prime}y_2 &= -H^x_{02}(\textbf{x}) \\
            y_1y_2^{\prime} + y_1^{\prime}y_2 &= -H^x_{12}(\textbf{x}) \\
            0&= G^x_2(\textbf{x}) y_1^2 + H^x_{12}(\textbf{x})y_1y_2 + G^x_1(\textbf{x})y_2^2 \\
            0&= G^x_0(\textbf{x}) y_1^2 + H^x_{01}(\textbf{x})y_0y_1 + G^x_1(\textbf{x})y_0^2 \\
            0&= G^x_0(\textbf{x}) y_2^2 + H^x_{02}(\textbf{x})y_0y_2 + G^x_2(\textbf{x})y_0^2
        \end{align*}
        Solving this system is merely a matter of case by case analysis, so we leave the interested reader to examine the PARI/gp script ordersearch.gp or see \cite[Section 6.2]{Hutz}.  Note that the method is valid only for $\text{char}(K) =p \neq 2$ due to using the quadratic formula for degree 2 polynomials, however, for $\mathbb{F}_2$ it is a simple matter to examine all possible points in $\mathbb{P}_x^2(\mathbb{F}_2) \times \mathbb{P}_y^2(\mathbb{F}_2)$.

    \section{Riemann Zeta Function of a Wehler K3 surface} \label{sect6}
        Let $p \in \mathbb{Z}$ be a prime and let $N_m = \#S(\mathbb{F}_{p^m})$ be the number of $\mathbb{F}_{p^m}$-rational points on $S$.  Using Section \ref{ssect2}, we can quickly determine the number of points in a fiber, and hence can efficiently determine $N_m$.  We can then determine the Riemann zeta function of a surface $S$ and give an example of a family of Wehler K3 surfaces that have Picard number $2$.

        \begin{defn}
            The \emph{Riemann zeta function} of $S$ is
            \begin{equation*}
                Z(S,T) = \exp\left(\sum_{m =1}^{\infty} N_m \frac{T^m}{m}\right),
            \end{equation*}
            where $\exp$ denotes exponentiation.
        \end{defn}

        The Riemann zeta function encodes much information about the set of points of $S$.  Knowing
        the Riemann zeta function therefore leads to knowing information about $S$.  In particular, we produce a K3 surface with Picard number two.  One way to
        compute the Riemann zeta function is to determine $N_m$ for enough $m$.  How many is
        enough depends on the variety.  Since $S$ is a K3 surface, it has Betti numbers
        $b_0=b_4=1$, $b_1=b_3=0$, and $b_2=22$, consequently, we will see that we need only the first 11 $N_m$s to determine the Riemann zeta function.

    \subsection{The Riemann Hypothesis}
        The Riemann zeta function satisfies the Riemann Hypothesis, as shown by
        Deligne \cite{Deligne2}.  We have that $\dim(S) =2$ and that $K = \mathbb{F}_{3}$.  We are using $p=3$ since the algorithm from Section \ref{ssect2} is not applicable in the case $p=2$.  Since $S$ is a K3 surface, we have that
        \begin{align}
            Z(S,T) &= \exp\left(\sum_{m =1}^{\infty} N_m \frac{T^m}{m}\right) =
            \frac{1}{P_0(T)P_2(T)P_{4}(T)} \notag\\
            &=\frac{1}{(1-T)(\prod_{i=1}^{22}(1-\alpha_iT))(1-9T)}\label{eq1}.
        \end{align}
        Our goal is to calculate the 22 $\alpha_i$s.  The method we use is to
        \begin{enumerate}
            \item Compute the $N_m$s programmatically.
            \item Write equation (\ref{eq1}) as a power series in $T$.
            \item Equate coefficients of the powers of $T$ in the expansion of equation (\ref{eq1}) with the $N_m$s.
        \end{enumerate}
        From the above information, we are able to calculate the power sums of the $\alpha_i$,
        from which we can compute the symmetric polynomials in the $\alpha_i$.  Consequently, we
        can determine the polynomial with the $\alpha_i$ as its roots.

        Fortunately, the 22 $\alpha_i$s are not all independent, and we show next that it is sufficient to determine only the first 11 $N_m$s.  The functional equation of the Riemann zeta function says that we can group the $\alpha_i$s into pairs which satisfy $\alpha_i\alpha_i^{\prime} = 9$.
        Hence, we get
        \begin{equation*}
            P_2(T) = \prod_{i=1}^{22}(1-\alpha_iT) = \prod_{i=1}^{11}(1-\alpha_iT)(1-\alpha^{\prime}_iT) = \prod_{i=1}^{11}(1-a_iT + 9T^2),
        \end{equation*}
        where
        \begin{equation*}
            a_i = \alpha_i + \frac{9}{\alpha_i} \quad \text{and} \quad \alpha_i\alpha^{\prime}_i=9.
        \end{equation*}
        Taking the natural log of both sides and expanding we get the
        following equations, where $P_j = \sum_{i=1}^{11}a_i^j$:
        \begin{align*}
            P_1&=N_1 - 10\\
            P_2&=N_2 +116\\
            P_3&=N_3 + 27P_1 - 730\\
            P_4&=N_4 + 36P_2 - 8344\\
            P_5&=N_5 + 45P_3 - 405P_1 - 59050\\
            P_6&=N_6 + 54P_4 - 729P_2 - 515404\\
            P_7&=N_7 + 63P_5 - 1134P_3 + 5103P_1 - 4782970\\
            P_8&=N_8 + 72P_6 - 1620P_4 + 11664P_2 - 43191064\\
            P_9&=N_9 + 81P_7 - 2187P_5 + 21870P_3 - 59049P_1 - 387420490\\
            P_{10}&=N_{10} + 90P_8 - 2835P_6 + 36450P_4 - 164025P_2 - 3485485324\\
            P_{11}&=N_{11} + 99P_9 - 3564P_7 + 56133P_5 - 360855P_3 + 649539P_1 - 31381059610
        \end{align*}
        From these 11 power sums, it is possible to recover the first 11 symmetric polynomials $C_1,\ldots,C_{11}$ in the $a_i$s using the Newton-Girard formulas \cite{Seroul}, which relates the $P_k$ to the $C_i$.  We now have the polynomial whose roots are the $a_i$, and we know that
        \begin{equation*}
            \alpha_i^2 - a_i \alpha_i + 9 =0.
        \end{equation*}
        Therefore, we can solve for
        \begin{equation*}
            \alpha_i = \frac{a_i \pm \sqrt{a_i^2 - 36}}{2}.
        \end{equation*}
        So we determine the $a_i$ and can then determine the $\alpha_i$. Therefore, we only need to determine the first 11 $N_m$s to determine $Z(S,T)$.
        \begin{rem}
            We could compute the $C_i$ for $p \neq 3$ if we could feasibly determine the $N_m$s in this case.
        \end{rem}

    \subsection{An Example Computation of the $N_m$s} \label{ssect3}
        Thanks to the computation cluster at Amherst College for providing CPU
        time, we determine the first 11 $N_m$s for the Wehler K3 surface
            \begin{align*}
               L&= [1,0,0,0,1,0,0,0,1]\\
               Q&=[2, 2, 2, 0, 1, 1, 1, 2, 1, 1, 2, 1, 0, 0, 2, 2, 2, 2, 0, 1, 0, 2, 0, 2, 2, 0, 2, 1, 1, 1, 0, 1, 0, 0, 0, 0]
            \end{align*}
            The $N_m$s are shown in Table \ref{table13}.
            \begin{table}
            \caption{Point Counting Data}
            \label{table13}
                \begin{center}
                \begin{tabular}{|c|c|}
                    \hline
                    $m$ & $N_m$ \\
                    \hline
                    1 & 13 \\
                    \hline
                    2 & 97 \\
                    \hline
                    3 & 784 \\
                    \hline
                    4 & 6877 \\
                    \hline
                    5 & 60238 \\
                    \hline
                    6 & 533440 \\
                    \hline
                    7 & 4782322\\
                    \hline
                    8 & 43047613\\
                    \hline
                    9 & 387464230 \\
                    \hline
                    10 & 3486563272\\
                    \hline
                    11 & 31382110828\\
                    \hline
                \end{tabular}
                \end{center}
            \end{table}
            From the $N_m$s we compute $P_2(T)$ to be
            \begin{align*}
                P_2(T)&=(9T^2 - 6T + 1)(3486784401T^{20} + 1162261467T^{19} + 258280326T^{18} + 43046721T^{17} \\
                &- 14348907T^{16} - 14348907T^{15} - 4782969T^{14} - 531441T^{13}+ 177147T^{12} + 118098T^{11}\\
                &+ 78732T^{10} + 13122T^9 + 2187T^8 - 729T^7 - 729T^6- 243T^5 - 27T^4 + 9T^3 + 6T^2 + 3T + 1).
            \end{align*}
            Examining the $\alpha_i$ we know from \cite[Corollary 2.3]{vanluijk} that the number of $\frac{\alpha_i}{3}$ which are roots of unity is an upper bound on the Picard number.  In this example, the upper bound is then 2.  However, we also know that the Picard number is at least 2, see \cite{Wehler}.  Hence, the surface has Picard number exactly 2.  Moreover, we have an infinite family of Wehler K3 surfaces with Picard number $2$, since the Riemann zeta function over $\mathbb{F}_3$ is unchanged by adding multiples of $3$ to the coefficients.

\providecommand\biburl[1]{\texttt{#1}}

\end{document}